\documentclass[a4paper, 12pt]{amsart}
\usepackage{amssymb, amsfonts, amsmath, bm}
\usepackage[all]{xy}
\usepackage[shortlabels]{enumitem}
\usepackage{hyperref, color}
\usepackage{commath}
\usepackage{esdiff}
\usepackage{graphicx,subfigure}
\usepackage[center]{caption}
\usepackage{amsthm}
\usepackage{diagbox}
\usepackage{multirow}
\usepackage{esdiff}
\usepackage[pagewise]{lineno}
\usepackage[utf8]{inputenc}	
\usepackage{hyperref, color}

\hypersetup{
	colorlinks = true,
	linkcolor = black,
	filecolor = black,
	urlcolor = black,
	citecolor = black
}

\newcommand{\Rbb}{\mathbb{R}}

\newcommand{\Nbb}{\mathbb{N}}

\newcommand{\tA}{\widetilde{A}}
\newcommand{\tH}{\widetilde{H}}

\newcommand{\bnabla}{\bar{\nabla}}
\newcommand{\bmu}{\bar{\mu}}

\newcommand{\bg}{\bar{g}}
\newcommand{\bK}{\bar{K}}
\newcommand{\bDelta}{\bar{\Delta}}
\newcommand{\bkappa}{\bar{\kappa}}

\newcommand{\bSigma}{\bar{\Sigma}}

\newcommand{\betaa}{\bar{\eta}}

\newcommand{\bx}{\bar{x}}

\DeclareMathOperator{\intt}{\mathrm{int}}
\DeclareMathOperator{\tr}{\mathrm{trace}}
\DeclareMathOperator{\divv}{\mathrm{div}}

\newtheorem{theorem}{Theorem}[section]
\newtheorem{lemma}[theorem]{Lemma}
\newtheorem{proposition}[theorem]{Proposition}
\newtheorem{corollary}[theorem]{Corollary}
\newtheorem{conjecture}[theorem]{Conjecture}

\theoremstyle{definition}

\newtheorem*{remark}{Remark}
\newtheorem*{acknowledgements}{Acknowledgements}

\numberwithin{equation}{section}
\numberwithin{figure}{section}

\begin{document}
	
	\title[A $\sigma$-homothetic uniqueness of the critical catenoid]{A $\sigma$-homothetic uniqueness of the critical catenoid}
	
	\author[Roney Santos, Iury Domingos and Feliciano Vit\'orio]{Roney Santos, Iury Domingos and Feliciano Vit\'orio}

    \address{Universidade Federal de Alagoas\\
            Av. Manoel Severino Barbosa S/N,
            57309-005 Arapiraca - AL,
            Brazil.}
	\email{iury.domingos@im.ufal.br}
	
	\address{Universidade de S\~ao Paulo\\
		Departamento de Matem\'atica\\
        Rua do Mat\~ao, 05508-900, S\~ao Paulo - SP, Brazil.}
	\email{roneypsantos@ime.usp.br}
	
	\address{Universidade Federal de Alagoas\\
		Instituto de Matem\'{a}tica\\
		Campus A. C. Sim\~{o}es, BR 104 - Norte, Km 97, 57072-970, Macei\'o - AL, Brazil.}	
	\email{feliciano@pos.mat.ufal.br}
	
	%
	
	\thanks{I. Domingos was partially supported by the Brazilian National Council for Scientific and Technological Development (CNPq), grant 409513/2023-7. R. Santos was supported by Instituto Serrapilheira, grant ``New perspectives of the min-max theory for the area functional''. F. Vit\'orio was partially supported by CNPq (405468/2021- 0)}
	
	
	\begin{abstract}
		We prove a uniqueness result for free boundary minimal annuli in the unit Euclidean three-ball that are $\sigma$-homothetic to the critical catenoid.
	\end{abstract}

    \maketitle
	
	\section{Introduction}
	
	
	A \textit{free boundary minimal immersion} into Euclidean $n$-dimensional ball $B^n$ centered at the origin of $\Rbb^n$ is a minimal isometric immersion $x: \Sigma \to B^n$ of a $k$-dimensional smooth manifold $\Sigma$ which meets the boundary of $B^n$ orthogonally.

    This kind of immersion has been receiving much attention due to the works of Fraser-Shoen \cite{fraser2011first, fraser2016sharp} about the Steklov eigenvalues on compact surfaces with non-empty boundary.
	
	
	
	The goal in these notes is to study the behavior of the conformal factor of two immersed free boundary minimal annuli in $B^3,$ and to conclude a uniqueness result in the conformal class of the critical catenoid. Here, given two metrics $g$ and $\bg$ in a smooth manifold $\Sigma$ such that $\bg = e^{2\varphi}g$ for some function $\varphi \in C^1(\Sigma),$ we are calling \textit{conformal factor} the function $\varphi.$ Our first result reads as follows.
	
	\begin{theorem}\label{maintheorem}
		Let $\Sigma$ and $\bSigma$ be two conformal immersed free boundary minimal annuli in $B^3$ with conformal factor $\varphi \in C^4({\Sigma}).$ Then, there exists $C \in \Rbb$ such that
		\[\begin{cases}
			\Delta \varphi =\left(1 - e^{C - 2\varphi}\right)K\ &\text{in}\ \Sigma\\
			\partial_\nu\varphi = e^{\varphi} - 1\ &\text{on}\ \partial\Sigma,
		\end{cases}\]
		where $\Delta$ and $K$ are respectively the Laplacian and the Gaussian curvature of $\Sigma.$
	\end{theorem}
	
	The Neumann boundary condition can be verified for conformal factors of any two conformal surfaces inside $B^n$ provided that both meet $\partial B^n$ orthogonally.

    Despite the celebrated Nitsche theorem \cite{nitsche1985stationary} classifying the flat equatorial disc as the unique free boundary minimal disc in $B^3,$ there are no known rigidity results for free boundary minimal surfaces in $B^3$ of other topological types, only under topological assumptions. In this sense, we want to prove a uniqueness result in the conformal class of the critical catenoid provided that the conformal factor is constant at least in one of the boundary components.

    \begin{theorem}\label{maintheorem2}
	An immersed free boundary minimal annulus in $B^3$ conformal to the critical catenoid whose conformal factor is constant along at least one of its boundary component and of class $C^4$ must be congruent to the critical catenoid.
    \end{theorem}
    
    Remember that two compact surfaces $\Sigma_1$ and $\Sigma_2$ with non-empty boundaries are $\sigma$-homothetic if there is a conformal diffeomorphism between them that is an homothety along the boundary. Hence, in particular, our result shows that a free boundary minimal annulus in $B^3$ that is $\sigma$-homothetic to the critical catenoid is indeed isometric to the critical catenoid.
    
    The concept of $\sigma$-homothety appears in \cite{fraser2012minimal, fraser2016sharp}, where the authors have studied the maximum of the first Steklov eigenvalue times the length of the boundary of a surface.

    To be more precise, let us consider an abstract surface $(\Sigma,g)$ with non-empty boundary. Let $\sigma_1(g)$ and $L(\partial\Sigma,g)$ be the first Steklov eigenvalue and the length of $\partial\Sigma$ both with respect to $g.$ When $\Sigma$ is annulus, Fraser-Schoen showed that the maximum of $\sigma_1(g)L(\partial\Sigma,g)$ over all smooth metric on $\Sigma$ is realised by surfaces that are $\sigma$-homothetic to the critical catenoid.
	
	There is an interesting analogy between closed minimal surfaces in the unit sphere $S^3$ and free boundary minimal surfaces in $B^3.$ This analogy have inspired several mathematicians, and during the last few years there has been a substantial development in the theory of free boundary minimal surfaces in $B^3.$ On this analogy, it is important to mention the acclaimed conjecture about the uniqueness of the critical catenoid due to Fraser-Li \cite{fraser2014compactnessof}. As observed by Li \cite{li2019free}, this conjecture was claimed by Nitsche \cite{nitsche1985stationary} without a proof.
	
	\begin{conjecture}\label{conjecture}
		An embedded free boundary minimal annulus in $B^3$ must be congruent to the critical catenoid.
	\end{conjecture}
	
	In comparison to the theory of closed minimal surfaces in $S^3,$ this conjecture is analogous to the uniqueness of the Clifford torus in $S^3,$ known as Lawson conjecture, and confirmed by Brendle \cite{brendle2013embedded}. The examples of \textit{self-intersecting} free boundary minimal annuli in $B^3$ constructed by Fernandez-Hauswirth-Mira \cite{fernandez2022free} and Kapouleas-McGrath \cite{kapouleas2022free} show that embeddedness is an essential hypothesis in the conjecture, as it was in the Lawson conjecture.
	
	Some advances on the validity of the conjecture was obtained, for instance, in the works of Ambrozio-Nunes \cite{ambrozio2016gap}, Devyver \cite{MR3928806}, Fraser-Schoen \cite{fraser2016sharp}, McGrath \cite{mcgrath2016characterization}, Kusner-McGrath \cite{kusner2020free}, Seo \cite{seo2021sufficient} and Tran \cite{MR4077171}. 
        

    Up to our knowledge, there is no result relating the conformal factor of two conformal free boundary minimal surfaces in $B^3$. Our Theorem \ref{maintheorem} shows that, at least for annuli, we can derive a relation of this kind. Also, in Theorem \ref{maintheorem2}, we are able to remark that this problem has a unique solution under an assumption on the boundary of the surface, which gives a geometric characterization.
 
    A free boundary minimal surface in $B^3$ is topologically an annulus if and only if it has no umbilical points. Because of this, except for the Neumann boundary condition, our proof of the results holds only for annuli. The assumption that one of the surfaces is the critical catenoid in Theorem \ref{maintheorem2} is necessary because, in this case, as we will prove, both surfaces have constant Gaussian curvatures along their boundaries.
	
	
	\begin{acknowledgements}
		We are grateful to J. Espinar for his suggestions about this work.
	\end{acknowledgements}

	\section{Preliminaries}
	
	For convenience of the reader, we collect here the results of our previous work \cite{domingos2021note} that we will use to study the conformal changes on free boundary minimal annuli in $B^3.$ They can be found in the more general context of free boundary hypersurfaces with constant mean curvature in balls of space forms in \cite{domingos2021note}.
	
	The first one is a relation between the extrinsic geometry along the boundary of a $n$-dimensional free boundary minimal hypersurface $\Sigma$ viewed from the original hypersurface and from the sphere that delimits $B^{n+1}$. We will denote by $A$ the shape operator of $\Sigma$ and by $\tA$ the shape operator of $\partial\Sigma$ as a hypersurface of $\partial B^{n+1}.$
	
	\begin{lemma}\label{dSigma}
		Let $\Sigma$ be an immersed free boundary minimal hypersurface in $B^{n+1}.$ Then, at the points of $\partial \Sigma$ we have
		\begin{itemize}
			\item[$\mathrm{(i)}$] $\displaystyle |A|^2 = |\widetilde{A}|^2 + \tH^2$;
			\item[$\mathrm{(ii)}$] $\partial_\nu|A|^2 = - 2 \left(|\tA|^2 + (n + 1)\tH^2\right),$
		\end{itemize}
		where $\nu$ is the exterior unit conormal to $\Sigma$ along $\partial \Sigma,$ and $\widetilde{H} = \tr(\tA).$
	\end{lemma}
	
	In particular, in the case $n=2,$ if $K$ denotes the Gaussian curvature of $\Sigma,$ then $\partial_\nu K = -4K.$
	
	\begin{remark}
		The second identity of Lemma \ref{dSigma} was first obtained in the work of Wheeler-Wheeler \cite{MR3874608}.
	\end{remark}
	
	The second result follows as a consequence of a formula that counts the umbilical points of a free boundary minimal surface in $B^3$ in terms of its Euler characteristic. For a proof, see, for instance, \cite[Lemma 4.3]{li2019free} or \cite[Corollary 4.4]{domingos2021note}.
	
	\begin{lemma}\label{umbiliccorollary}
		An immersed free boundary minimal surface in $B^3$ is an annulus if and only if it has no umbilical points.
	\end{lemma}
	
	An important example of free boundary minimal surface in $B^3$ is the piece of catenoid contained in $B^3$ that intersects $\partial B^3$ orthogonally called the \textit{critical catenoid}. We can parametrise it by the immersion $x: [-T_0, T_0] \times [0, 2\pi] \to B^3$ given by
	\[x(t, \theta) = a(\cosh t \cos\theta, \cosh t \sin\theta, t),\]
	where $T_0$ is the unique positive solution of $t\tanh t = 1$ and the dilatation constant is $a = (T_0\cosh T_0)^{-1}.$ By means of straightforward computations, one infers that the Gaussian curvature $K$ of the critical catenoid is given by
	\[K = - \frac{1}{a^2 \cosh^4t}.\]
	
	It is a well-known result that the only immersed free boundary minimal surface in $B^3$ whose at least one of its boundary curves is a circle must be congruent either to the equatorial disk or to the critical catenoid (see, for instance, Pyo \cite{pyo2010minimal}). Moreover, if $\Sigma$ is a free boundary surface in $B^3,$ then its boundary components are line of curvatures of both $\Sigma$ and $\partial B^3$ by the Joachimsthal theorem. The following result is a consequence of these two facts because the boundary curves of $\Sigma$ are contained in the sphere.
	
	\begin{lemma}\label{constantcurvature}
		An immersed free boundary minimal surface in $B^3$ whose Gaussian curvature is constant along at least one of its boundary curves is congruent either to the equatorial disk or to the critical catenoid.
	\end{lemma}
	
	\section{Conformal changes and free boundary condition}\label{conformalchanges}
	We begin this section by recalling that the condition that a surface meets $\partial B^n$ orthogonally imposes a strong restriction on the geometry of the boundary of the surface.

    \begin{proposition}\label{geodesiccurvature}
        Let $x:(\Sigma,g) \to B^n$ be an isometrically immersed surface into $B^n$ meeting $\partial B^n$ orthogonally. Then, the geodesic curvature of the boundary components of $\Sigma$ in the exterior direction are equals $-1.$
    \end{proposition}

\begin{proof}
    Let $\nu$ denotes the unit exterior conormal of $\Sigma$ along $\partial\Sigma.$ Since $dx(\nu) = x$ by the meeting condition along $\partial\Sigma,$ we apply Gauss equation to conclude that
        \[(x^*\nabla)_X\nu = x^*\left(\nabla_{dx(X)} x\right) = X\]
        for each vector field $X$ tangent to $\partial\Sigma,$ where $\nabla$ and $x^*\nabla$ denote respectively the connection of $x(\Sigma)$ and the pullback connection of $\nabla.$ Thus, the geodesic curvature $\kappa$ of $\partial \Sigma$ in the direction of $\nu$ is
        \[\kappa = -g(\nabla_T\nu, T) = -1,\]
        where $T$ is the unit tangent of a boundary component of $\Sigma.$
\end{proof}
    
	In particular, the boundary geometry implies that the conformal factor of two conformal surfaces in $B^n$ meeting $\partial B^n$ orthogonally is such that it must respect the following equation at the boundary points.
	
	\begin{proposition}\label{normalderivative}
		Let $x:(\Sigma,g) \to B^n$ and $\bx:(\Sigma,\bg) \to B^n$ be two isometrically immersed surfaces into $B^n$ meeting $\partial B^n$ orthogonally and such that $\bg = e^{2\varphi}g$ for some $\varphi \in C^1(\Sigma).$ Then, at the points of $\partial \Sigma$ we have
		\[\partial_\nu\varphi = e^\varphi - 1,\]
		where $\nu$ denotes the unit exterior conormal to $\Sigma$ along $\partial \Sigma$ with respect to $g.$
	\end{proposition}
	
	\begin{proof}
		Let $\kappa$ and $\bkappa$ denote the geodesic curvatures of $\partial\Sigma$ respectively as curves inside $(\Sigma,g)$ and $(\Sigma,\bg)$ in the exterior direction. By Proposition \ref{geodesiccurvature}, we have $\kappa = \bkappa = -1.$ Replacing this in the formula
		\[e^\varphi\bkappa = \kappa - \partial_\nu\varphi\]
		we obtain $e^\varphi = 1 + \partial_\nu\varphi,$ as claimed.
	\end{proof}
	
	\begin{remark}
		The formula of Proposition \ref{normalderivative} can be used to prove the inequality of Fraser-Schoen \cite{fraser2011first, fraser2012minimal} about the length of the boundary of a free boundary minimal surface in $B^n$ in its conformal orbit. More precisely, we can show that if $\Sigma$ is an immersed free boundary minimal surface in $B^n$ and $f: B^n \to B^n$ is a conformal diffeomorphism, then
		\[L(\partial\Sigma) \geq L(f(\partial\Sigma)).\]
		
		To show this, let $\Phi \in C^\infty(B^n)$ be the conformal factor of $f.$ Thus, there exists $x_0 \in \Rbb^n$ with $|x_0| > 1$ such that
		\[\Phi(x) = \log \frac{|x_0|^2 - 1}{|x - x_0|^2}.\]
		If $\nabla$ and $\Delta$ are the gradient and the Laplacian of $\Sigma$ respectively, one can check that
		\[\nabla|x-x_0|^2 = 2(x-x_0)^T\ \ \ \mbox{and}\ \ \ \Delta|x-x_0|^2 = 4,\]
		where $(\cdot)^T$ stands for the orthogonal projection onto the tangent space of $\partial \Sigma,$ and the second equality is consequence of the minimality of $\Sigma.$ Therefore, setting $\varphi = \Phi|_\Sigma,$ we have
		\[\Delta \varphi = \frac{4|(x-x_0)^T|^2}{|x-x_0|^4} - \frac{4}{|x-x_0|^2} \leq 0.\]
		Combining divergence theorem with Proposition \ref{normalderivative}, we obtain that
		\[0 \geq \int_\Sigma \Delta\varphi\, d\mu = \int_{\partial\Sigma} \partial_\nu \varphi\, ds = L(f(\partial\Sigma)) - L(\partial\Sigma).\]

		In \cite{fraser2012minimal}, the authors choose a vector field $V$ given by
		\[V = \frac{x - x_0}{|x - x_0|^2}.\]
		This is not necessarily a vector field tangent to $\Sigma$ at each point. But considering the divergence operator $\divv_\Sigma V = \langle \bnabla_{v_1} V, v_1\rangle + \langle \bnabla_{v_2} V, v_2\rangle,$ where $\{v_1, v_2\}$ is an orthonormal basis while $\bnabla$ and $\langle \cdot, \cdot\rangle$ are respectively the connection and the metric of $\Rbb^n,$ one can see by means of a direct computation that $\Delta\varphi = -2\divv_\Sigma V.$
	\end{remark}
	
	The next result gives the behavior of the conformal factor between two free boundary minimal annuli in $B^3.$ The conformal factor must satisfies a differential equation involving the Gaussian curvature of one of the surfaces.
	
	\begin{theorem}\label{laplacianphi}
		Let $x: (\Sigma,g) \to B^3$ and $\bx:(\Sigma,\bg) \to B^3$ be two free boundary minimal immersions of annuli into $B^3$ such that $\bg = e^{2\varphi}g$ for some $\varphi \in C^4({\Sigma}).$ Then, there exists $C \in \Rbb$ such that
        \[\begin{cases}
			\Delta \varphi =\left(1 - e^{C - 2\varphi}\right)K\ &\text{in}\ \Sigma\\
			\partial_\nu\varphi = e^{\varphi} - 1\ &\text{on}\ \partial\Sigma,
		\end{cases}\]
		where $K$ is the Gaussian curvature of $\Sigma.$        
	\end{theorem}
	
	\begin{proof}
		By conformality between $g$ and $\bg,$ we get $e^{2\varphi} \bK = K - \Delta\varphi,$ where $\bK$ is the Gaussian curvature of $\bg.$ Moreover, by Lemma \ref{umbiliccorollary}, we have that $K$ and $\bK$ does not vanish at any point of $\Sigma.$ Hence, $4K = \Delta \log(-K)$ and $4\bK = \bDelta\log(-\bK)$ in $\Sigma$ by minimality, where $\bDelta$ is the Laplacian of $\bg.$ Therefore, one can check that
		\[4\Delta \varphi = \Delta \log \left(\bK^{-1} K\right)\]
		at each point of $\Sigma.$ We consider the function $f\in C^2(\Sigma)$ given by $f = 4\varphi - \log\left(\bK^{-1}K\right).$ It follows by Lemma \ref{dSigma} that
		\[\partial_\nu\log\left(\bK^{-1}K\right) = K^{-1}\partial_\nu K - \bK^{-1}\partial_\nu \bK = 4(e^\varphi-1).\]
		Thus, by Proposition \ref{normalderivative}, we have $\partial_\nu f = 0$ along $\partial \Sigma.$ Since $f$ is harmonic, then $f$ is constant, which means that 
		\[4\varphi = C + \log(\bK^{-1} K)\]
		for some $C \in \Rbb,$ and this prove the first part of our assertion. By means of a direct computation, we obtain $\Delta \varphi =\left(1 - e^{C - 2\varphi}\right)K.$
        \end{proof}

        Our main result is a direct consequence of the proof of the previous result.

 \begin{corollary}
     Let $x: (\Sigma,g) \to B^3$ and $\bx:(\Sigma,\bg) \to B^3$ be two free boundary minimal immersions of annuli into $B^3$ such that $\bg = e^{2\varphi}g$ for some $\varphi \in C^4({\Sigma}).$ Suppose $x(\Sigma)$ is the critical catenoid and $\varphi$ is constant along at least one boundary component of $\Sigma.$ Then, $\bx(\Sigma)$ is congruent to the critical catenoid.
 \end{corollary}
 
     \begin{proof}
        We know that $4\varphi = C + \log(\bK^{-1}K)$ for some real constant $C.$ Since $\Sigma$ is the critical catenoid and $\varphi$ is constant along a boundary component $\Gamma$ of $\Sigma$ by hypothesis, we must have $\bK$ constant along $\Gamma.$ Hence, by Lemma \ref{constantcurvature} $\bx(\Sigma)$ is isometric to the critical catenoid.
    \end{proof}
	
	\bibliographystyle{siam}
	\bibliography{references}

\begin{thebibliography}{10}

\bibitem{ambrozio2016gap}
{\sc L.~Ambrozio and I.~Nunes}, {\em A gap theorem for free boundary minimal
  surfaces in the three-ball}, Communications in Analysis and Geometry, 29
  (2016).

\bibitem{brendle2013embedded}
{\sc S.~Brendle}, {\em Embedded minimal tori in {$S^3$} and the {L}awson
  conjecture}, Acta Math., 211 (2013), pp.~177--190.

\bibitem{MR3928806}
{\sc B.~Devyver}, {\em Index of the critical catenoid}, Geom. Dedicata, 199
  (2019), pp.~355--371.

\bibitem{domingos2021note}
{\sc I.~Domingos, R.~Santos, and F.~Vit\'{o}rio}, {\em A note on free boundary
  hypersurfaces in space form balls}, Arch. Math. (Basel), 121 (2023),
  pp.~197--209.

\bibitem{fernandez2022free}
{\sc I.~Fern\'{a}ndez, L.~Hauswirth, and P.~Mira}, {\em Free boundary minimal
  annuli immersed in the unit ball}, Arch. Ration. Mech. Anal., 247 (2023),
  pp.~Paper No. 108, 44.

\bibitem{fraser2014compactnessof}
{\sc A.~Fraser and M.~M.-c. Li}, {\em Compactness of the space of embedded
  minimal surfaces with free boundary in three-manifolds with nonnegative
  {R}icci curvature and convex boundary}, J. Differential Geom., 96 (2014),
  pp.~183--200.

\bibitem{fraser2011first}
{\sc A.~Fraser and R.~Schoen}, {\em The first {S}teklov eigenvalue, conformal
  geometry, and minimal surfaces}, Adv. Math., 226 (2011), pp.~4011--4030.

\bibitem{fraser2012minimal}
\leavevmode\vrule height 2pt depth -1.6pt width 23pt, {\em Minimal surfaces and
  eigenvalue problems}, in Geometric analysis, mathematical relativity, and
  nonlinear partial differential equations, vol.~599 of Contemp. Math., Amer.
  Math. Soc., Providence, RI, 2013, pp.~105--121.

\bibitem{fraser2016sharp}
\leavevmode\vrule height 2pt depth -1.6pt width 23pt, {\em Sharp eigenvalue
  bounds and minimal surfaces in the ball}, Invent. Math., 203 (2016),
  pp.~823--890.

\bibitem{kapouleas2022free}
{\sc N.~Kapouleas and P.~McGrath}, {\em Free boundary minimal annuli immersed
  in the unit 3-ball}, arXiv preprint arXiv:2212.09680,  (2022).

\bibitem{kusner2020free}
{\sc R.~Kusner and P.~McGrath}, {\em On free boundary minimal annuli embedded
  in the unit ball}, arXiv preprint arXiv:2011.06884,  (2020).

\bibitem{li2019free}
{\sc M.~M.-c. Li}, {\em Free boundary minimal surfaces in the unit ball: recent
  advances and open questions}, in Proceedings of the {I}nternational
  {C}onsortium of {C}hinese {M}athematicians 2017, Int. Press, Boston, MA,
  [2020] \copyright 2020, pp.~401--435.

\bibitem{mcgrath2016characterization}
{\sc P.~McGrath}, {\em A characterization of the critical catenoid}, Indiana
  Univ. Math. J., 67 (2018), pp.~889--897.

\bibitem{nitsche1985stationary}
{\sc J.~C.~C. Nitsche}, {\em Stationary partitioning of convex bodies}, Arch.
  Rational Mech. Anal., 89 (1985), pp.~1--19.

\bibitem{pyo2010minimal}
{\sc J.~Pyo}, {\em Minimal annuli with constant contact angle along the planar
  boundaries}, Geometriae Dedicata, 146 (2010), pp.~159--164.

\bibitem{seo2021sufficient}
{\sc D.~Seo}, {\em Sufficient symmetry conditions for free boundary minimal
  annuli to be the critical catenoid}, arXiv preprint arXiv:2112.11877,
  (2021).

\bibitem{MR4077171}
{\sc H.~Tran}, {\em Index characterization for free boundary minimal surfaces},
  Comm. Anal. Geom., 28 (2020), pp.~189--222.

\bibitem{MR3874608}
{\sc G.~Wheeler and V.-M. Wheeler}, {\em Minimal hypersurfaces in the ball with
  free boundary}, Differential Geom. Appl., 62 (2019), pp.~120--127.

\end{thebibliography}
	
\end{document}